\documentclass{amsart}
\usepackage{amsmath,amsthm,amssymb,amsfonts}
\usepackage{graphicx} 
\usepackage{epstopdf}

\makeatletter

\def\proof{\@ifstar{P\,r\,o\,o\,f}{P\,r\,o\,o\,f.\ }}

\renewcommand\th@remark{%
  \thm@headfont{\bfseries}%
  \normalfont 
  \thm@preskip\topsep \divide\thm@preskip\tw@
  \thm@postskip\thm@preskip
}

\renewenvironment{equation}{\refstepcounter{equation}$$}{\eqno{(\thesection.\theequation)}$$}

\makeatother

\newcounter{Example}[section]
\newcounter{Th}[section] \newcounter{Pr}[section] \newcounter{Lm}[section]
\newcounter{Remark}[section]
\newcounter{Def}[section]
\newcounter{lcounter}[section]
\newcounter{Corol}[section]

\newenvironment{Th}[1][\relax]
    {\medspace\refstepcounter{Th}T\,h\,e\,o\,r\,e\,m \arabic{section}.\theTh.\ \it}
    {\rm\medspace}

\newenvironment{Th.}[1][\relax]
    {\medspace\refstepcounter{Th}T\,h\,e\,o\,r\,e\,m \arabic{section}.\theTh.\ \it}
    {\rm\medspace}

\newenvironment{Pr.}[1][\relax]
    {\medspace\refstepcounter{Pr}P\,r\,o\,p\,o\,s\,i\,t\,i\,o\,n \arabic{section}.\thePr.\ \it}
    {\rm\medspace}

\newenvironment{Lm}[1][\relax]
    {\medspace\refstepcounter{Lm}L\,e\,m\,m\,a \arabic{section}.\theLm.\ \it}
    {\rm\medspace}

\newenvironment{Corol}[1][\relax]
    {\medspace\refstepcounter{Corol} C\,o\,r\,o\,l\,l\,a\,r\,y \arabic{section}.\theCorol.\ \it}
    {\rm\medspace}

\newenvironment{Def}[1][\relax]
    {\medspace\refstepcounter{Def}D\,e\,f\,i\,n\,i\,t\,i\,o\,n \arabic{section}.\theDef.\rm\ }
    {\medspace}

\newenvironment{Def.}[1][\relax]
    {\medspace\refstepcounter{Def}D\,e\,f\,i\,n\,i\,t\,i\,o\,n \arabic{section}.\theDef.\rm\ }
    {\medspace}

\def\au#1{\emph{#1}}
\def\tit#1{{#1}}
\def\R{{\mathbb R}}  
\def\C{{\mathbb G}}  
\def\cl{\mathop{\rm cl\,}}

\def\co{\mbox{\rm co}\,}

\def\aff {\mbox{\rm aff}\,}

\def\d {\partial\,}
\def\ep{\varepsilon} 
\def\Int {\mbox{\rm int\,}} 

\def\L {{\mathcal L}}  
\def\diam {\mbox{\rm diam}\,}
\def\diff {\,\frac{\,*\,}{}\,}

\begin{document}

\title[Polyhedral  approximations of strictly convex compacta]{Polyhedral approximations of strictly convex compacta}

\author{Maxim V. Balashov and Du\v{s}an Repov\v{s}}


\address{Department of Higher Mathematics, Moscow Institute of Physics and Technology, Institutski str. 9,
Dolgoprudny, Moscow region, Russia 141700. balashov@mail.mipt.ru}
\address{Faculty of Mathematics and Physics, and Faculty of Education, University of Ljubljana, Jadranska 19, Ljubljana, Slovenia 1000.
dusan.repovs@guest.arnes.si}



\begin{abstract}
We consider polyhedral approximations of strictly convex compacta
in finite dimensional Euclidean spaces (such compacta are also
uniformly convex). We obtain the best possible estimates for
errors of considered approximations in the Hausdorff metric. We
also obtain new estimates of an approximate algorithm for finding
the convex hulls.

\end{abstract}

\date{\today}
\maketitle

\section{Introduction}

\def\i {\mbox{\rm int}\,}

We begin by some definitions for a finite dimensional Euclidean
space $(\R^{n}, \|\cdot\|)$ over $\R$ with an inner product
$(\cdot,\cdot)$. Let $B_{r}(a)=\{ x\in \R^{n}\ |\ \| x-a\| \le
r\}$. Let $\cl A$ denote the \it closure \rm and $\Int A$ the \it
interior \rm of the subset $A\subset \R^{n}$. The \it diameter
\rm of the subset $A\subset \R^{n}$ is defined as $\diam A =
\sup\limits_{x,y\in A} \| x-y\|$. The \it distance \rm from the
point $x\in \R^{n}$ to the set $A\subset\R^{n}$ is given by the
formula $\varrho (x,A)=\inf\limits_{a\in A}\| x-a\|$. We shall
denote the \it convex hull \rm of a set $A\subset \R^{n}$ by $\co
A$, the convex hull of a function $f:\R^{n}\to \R$ by $\co f$ (cf.
\cite{Aubin,Polovinkin+Balashov,Rockafellar}).

The \it Hausdorff distance \rm between two subsets $A,B\subset
\R^{n}$ is defined as follows $h(A,B)=$
$$ = \max \left\{ \sup_{a\in A}\ \inf_{b\in B} \|a-b\| ,
\quad \sup_{b\in B}\ \inf_{a\in A} \|a-b\| \right\}=\inf\{ r>0\
|\ A\subset B+B_{r}(0),\ B\subset A+B_{r}(0) \}.
$$

The \it supporting function \rm of the subset $A\subset \R^{n}$
is defined as follows
\begin{equation}\label{**1}
s(p,A) = \sup\limits_{x\in A}(p,x),\qquad\forall p\in \R^{n}.
\end{equation}
The supporting function of any set $A$ is always lower
semicontinuous, positively uniform and convex. If the set $A$ is
bounded then the supporting function is Lipschitz continuous
\cite{Aubin,Polovinkin+Balashov}.

It follows from the separation theorem that for any convex
compacta $A,B$ in $\R^{n}$ (cf. \cite[Lemma
1.11.4]{Polovinkin+Balashov})
\begin{equation}\label{*3}
h(A,B) = \sup\limits_{\| p\|=1}|s(p,A)-s(p,B)|.
\end{equation}

A convex compactum in $\R^{n}$ is called \it strictly \rm convex
if its boundary contains no nontrivial line segments.

\begin{Def}\label{modulus} (Polyak \cite{Polyak}).
Let $E$ be a Banach space and let a subset $A\subset E$ be convex
and closed. \it The modulus of convexity \rm $\delta_{A}:\
[0,\diam A)\to [0,+\infty) $ is the function defined by
$$
\delta_{A}(\ep) = \sup\left\{ \delta\ge 0\ \left|\
B_{\delta}\left( \frac{x_{1}+x_{2}}{2}\right)\right.\subset A,\
\forall x_{1},x_{2}\in A:\ \| x_{1}-x_{2}\|=\ep \right\}.
$$\rm
\end{Def}

\begin{Def}\label{RM} (Polyak \cite{Polyak}).
Let $E$ be a Banach space and let a subset $A\subset E$ be convex
and closed. If the modulus of convexity $\delta_{A}(\ep )$ is
strictly positive for all $\ep\in (0,\diam A)$, then we call the
set $A$ \it uniformly convex (with modulus
$\delta_{A}(\cdot)$).\rm
\end{Def}

We proved in \cite{Balashov+Repovs2} that every uniformly convex
set is bounded and if the Banach space  $E$ contains a
nonsingleton uniformly convex set then it admits a uniformly
convex equivalent norm. We also proved that the function $\ep\to
\delta_{A}(\ep)/\ep$ is increasing (see also \cite[Lemma
1.e.8]{Lindestrauss+tzafriri}), and for any uniformly convex set
$A$ there exists a constant $C>0$ such that $\delta_{A}(\ep)\le
C\ep^{2}$.

The class of strictly convex compacta coincides with the class of
uniformly convex compacta with moduli of convexity
$\delta_{A}(\ep)>0$ for all permissible $\ep>0$ in the finite
dimensional case (cf. \cite{Balashov+Repovs2}).

\begin{Def}\label{grid} (\cite{sca, Polovinkin+Balashov}).
A \it grid \rm $\C$ with \it step \rm $\Delta\in (0,\frac12 )$ is
a finite collection of unit vectors $\{ p_{i}\}\subset \R^{n}$,
$i\in \overline{1,I}=\{1,\dots,I\}$, such that for any vector
$p\ne 0$, $p\in\R^{n}$, with $\frac{p}{\| p\|}\notin \C$ there
exist a set of indexes $I_{p}\subset \overline{1,I}$ and numbers
$\alpha_{i}>0$, $i\in I_{p}$, with the property
\begin{equation}\label{*2}
p=\sum\limits_{i\in I_{p}}\alpha_{i}p_{i},\qquad p_{i}\in\C,
\end{equation}
\begin{equation}\label{*1}
\| p_{i}-p_{j}\|<\Delta, \qquad\forall i,j\in I_{p}.
\end{equation}
\end{Def}

It is well known \cite{Aubin,Polovinkin+Balashov,Rockafellar}
that for any convex closed subset $A\subset \R^{n}$ we have
$$
A = \{ x\in\R^{n}\ |\ (p,x)\le s(p,A),\quad \forall p\in \d
B_{1}(0)\}.
$$

We shall consider external polyhedral approximation of the compact
$A\subset \R^{n}$ on the grid $\C$ from Definition 1.\ref{grid}
$$
\hat A = \{ x\in \R^{n}\ |\ (p,x)\le s(p,A),\quad \forall
p\in\C\}.
$$

From the inclusion $\C\subset \d B_{1}(0)$ we easily see that
$A\subset \hat A$. For an arbitrary convex compact set
$A\subset\R^{n}$ we have $h(A,\hat A)\le 2h(\{0\},A)\Delta$ (cf.
\cite{sca,Polovinkin+Balashov}). If $A=\bigcap\limits_{x\in
X}B_{R}(x)\ne\emptyset$, then $h(A,\hat A)\le 2R\Delta^{2}$ (cf.
\cite{sca,Polovinkin+Balashov}). Further we shall consider the
approximation of an arbitrary strict\-ly=uni\-form\-ly convex
compact set $A\subset\R^{n}$. Our further goals are\\
{\bf 1)} Estimate the error $h(A,\hat A)$ via the geometric
properties of the set $A$;\\
{\bf 2)} Suppose that we know a presupporting function $f(p)$ of
the convex compact $A$, i.e. the function $f$ is positively
uniform, continuous and $\co f(p) = s(p,A)$. Let
$$
\tilde A = \{ x\in\R^{n}\ |\ (p,x)\le f(p),\quad \forall p\in\C\}.
$$
Estimate the error $h(A,\tilde A)$ via the properties of the
function $f$ and geometric properties of the set $A$. In this
case we do not know the supporting function of the set $A$, but we
can find information about some properties of the set $A$:
diameter, modulus of convexity,
etc.\\
{\bf 3)} We shall consider an algorithm for calculating the convex
hull of a positively uniform function defined on the grid and
discuss estimates for the errors of such algorithms.

\section{Approximation by supporting functions}

\begin{Lm}\label{L1}  For a given grid $\C$
(Definition 1.\ref{grid}) with step $\Delta\in (0,\frac12 )$ in
the representation of any vector $p\ne 0$, $\frac{p}{\|
p\|}\notin \C$ by formulae (1.\ref{*2}), (1.\ref{*1}) the
following estimates hold
$$
\| \hat p - p_{j}\|<\Delta,\quad \forall j\in I_{p},\qquad 1\ge \|
\hat p \|\ge 1-\frac12\Delta^{2},
$$
where
\begin{equation}\label{*4}
\hat p=\frac{p}{\alpha}=\sum\limits_{i\in
I_{p}}\hat\alpha_{i}p_{i},\qquad \alpha =\sum\limits_{i\in
I_{p}}\alpha_{i},\quad \hat\alpha_{i}=\frac{\alpha_{i}}{\alpha}.
\end{equation}
\end{Lm}

\proof  From the definitions of $\hat p$ and $\hat \alpha_{i}$ we
obtain that
$$
\hat p = \sum\limits_{i\in I_{p}}\hat\alpha_{i}p_{i},\quad
\sum\limits_{i\in I_{p}}\hat\alpha_{i}=1.
$$
Hence
$$
\sum\limits_{i\in I_{p}}\hat\alpha_{i}(\hat p - p_{i})=0.
$$
By the triangle inequality we get
$$
\| \hat p - p_{j}\|\le \sum\limits_{i\in I_{p}}\hat\alpha_{i}\|
p_{i} -  p_{j}\|<\Delta,\qquad \forall j\in I_{p},
$$
\begin{equation}\label{*5}
\| \hat p\|\le \sum\limits_{i\in I_{p}}\hat\alpha_{i}\| p_{i}\|=1.
\end{equation}

The condition $\|p_i-p_j\|<\Delta$ is equivalent to the condition
$(p_i,p_j)=\frac12(\|p_i\|^2+\|p_j\|^2-\|p_i-p_j\|^2)\ge
1-\Delta^2/2$. Thus
 $1\ge \| \hat p\|\ge \|\hat p\|^2=\sum_{i,j}\hat\alpha_i\hat\alpha_j(p_i,p_j)\geq
(1-\frac12\Delta^2)\sum_{i,j}\hat\alpha_i\hat\alpha_j=1-\frac12\Delta^2$.
\qed

\def\CC {\mathcal C}
\def\UU {\mathcal U}

\begin{Def}\label{gridOp}
Let $f:\R^{n}\to\R$ be a positively uniform function. Let $\C$ be
a grid with step $\Delta\in (0,\frac12 )$. Define the grid
operators
$$
\CC f(p)=\left\{\begin{array}{lc} f(p),\quad \frac{p}{\|
p\|}\in\C,\\
\ \\+\infty,\quad \frac{p}{\| p\|}\notin\C,
\end{array}\right.\qquad\quad
 \UU f(p)=\left\{\begin{array}{lc} f(p),\quad \frac{p}{\|
p\|}\in\C,\\
\ \\ \sum\limits_{i\in I_{p}}\alpha_{i}f(p_{i}),\quad \frac{p}{\|
p\|}\notin\C.
\end{array}\right.
$$
Indices $I_{p}$ and numbers $\alpha_{i}$ are from Definition
1.\ref{grid}.
\end{Def}

\begin{Lm}\label{L2} (\cite[Lemma 5]{sca},
\cite[Lemma 2.6.2]{Polovinkin+Balashov}). Let $f:\R^{n}\to\R$ be a
positively uniform function.\\
(1) If the function $f$ is convex then $\CC f\ge f$,
$\co \CC f(p)=f(p)$ $\forall p\in\C$.\\
(2) If the function $f$ is convex then $f\le\co \UU f$.\\
(3) $\co \CC f = \co\UU f $, $\forall f$.\\
(4) $\co f\le \co \CC f$, $\forall f$.
\end{Lm}



The next lemma is a modification of Lemma 2.2 from
\cite{Balashov+Repovs2}.

\begin{Lm}\label{L3}
 Let  $A\subset \R^{n}$ be compact and uniformly convex set with the
 modulus of convexity $\delta$. Let $\ep\in (0,\diam A)$,
$\Delta\in (0,\frac12 )$. Let $p_{1},p_{2}\in \R^{n}$, $\|
p_{1}\|=1$, $1-\frac12\Delta^{2}\le \| p_{2}\|\le 1$. Let
$x_{i}=\arg\max\limits_{x\in A}(p_{i},x)$, $i=1,2$. If $\|
p_{1}-p_{2}\|<(4-\Delta^{2})\frac{\delta (\ep)}{\ep}$, then $\|
x_{1}-x_{2}\|<\ep$.
\end{Lm}

\proof Suppose that $\| x_{1}-x_{2}\|\ge \ep$. Let $t = \delta (\|
x_{1}-x_{2}\|)\ge \delta (\ep)$. By the condition
$$
B_{t}\left( \frac{x_{1}+x_{2}}{2}\right)\subset A
$$
we have that $(p_{1},x_{1})\ge
\left(p_{1},\frac{x_{1}+x_{2}}{2}\right)+t$,
\begin{equation}\label{*6}
(p_{1},x_{1}-x_{2})\ge 2t,
\end{equation}
$$
(p_{2},x_{2})\ge\left(p_{2},\frac{x_{1}+x_{2}}{2}\right)+t\| p_{2}
\|\ge
\left(p_{2},\frac{x_{1}+x_{2}}{2}\right)+t\left(1-\frac12\Delta^{2}\right),
$$
\begin{equation}\label{*7}
(p_{2},x_{2}-x_{1})\ge \left(2-\Delta^{2}\right)t.
\end{equation}
By formulae (2.\ref{*6}), (2.\ref{*7}) we obtain that
$$
(p_{1}-p_{2},x_{1}-x_{2})\ge (4-\Delta^{2})t,
$$
and
$$
\| p_{1}-p_{2}\|\ge (4-\Delta^{2})\frac{\delta (\|
x_{1}-x_{2}\|)}{\| x_{1}-x_{2}\|}.
$$
By Lemma 2.1 of \cite{Balashov+Repovs2} we have the inequality
$\frac{\delta (\| x_{1}-x_{2}\|)}{\| x_{1}-x_{2}\|}\ge
\frac{\delta (\ep)}{\ep}$ and
$$
\| p_{1}-p_{2}\|\ge (4-\Delta^{2})\frac{\delta (\ep)}{\ep}.
$$
Contradiction.\qed

\begin{Corol}\label{S1}
Let under the conditions of Lemma 2.\ref{L3} $\ep (\Delta)$ be a
solution of the equation $\frac{\delta
(\ep)}{\ep}=\frac{\Delta}{4-\Delta^{2}}$. If $\|
p_{1}-p_{2}\|<\Delta$ then $\| x_{1}-x_{2}\|<\ep (\Delta)$.
\end{Corol}

\proof The proof follows from Lemma
 2.\ref{L3} and strict monotonicity of the function $\frac{\delta
(\ep)}{\ep}$ \cite[Lemma 1.2]{Balashov+Repovs3}.\qed

\begin{Th}\label{T1}
 Let $p_{1},p_{2}\in \R^{n}$, $\|
p_{1}\|=1$, $1-\frac12\Delta^{2}\le\| p_{2}\|\le 1$, $\|
p_{1}-p_{2}\|<\Delta$. Let $A\subset \R^{n}$ be compact and
uniformly convex set with modulus of convexity $\delta$ and
$\Delta\in (0,\frac12 )$, $\delta (\diam A)/\diam
A>\frac{\Delta}{4-\Delta^{2}}$. Let $x_{i}=\arg\max\limits_{x\in
A}(p_{i},x)$, $i=1,2$. Then
\begin{equation}\label{*8}
s(p_{1},A)-s(p_{2},A)=(x_{2},p_{1}-p_{2})+\ep_{1}(\|
p_{1}-p_{2}\|)\| p_{1}-p_{2}\|,
\end{equation}
\begin{equation}\label{*9}
s(p_{2},A)-s(p_{1},A)=(x_{1},p_{2}-p_{1})+\ep_{2}(\|
p_{1}-p_{2}\|)\| p_{1}-p_{2}\|,
\end{equation}
and
$$
\max\{ |\ep_{1}(\| p_{1}-p_{2}\|)|\| p_{1}-p_{2}\|,\ |\ep_{1}(\|
p_{1}-p_{2}\|)|\| p_{1}-p_{2}\|\}\le \ep (\Delta)\Delta,
$$
where $\ep (\Delta)$ is  a solution of the equation $\frac{\delta
(\ep)}{\ep}=\frac{\Delta}{4-\Delta^{2}}$.
\end{Th}

\proof Equations (2.\ref{*8}) and  (2.\ref{*9}) are equivalent to
the condition of continuous gradient for (convex) supporting
function at the points $p_{2}$ and $p_{1}$, respectively. It is a
well known fact that the supporting function of the strictly
convex compact is continuously differentiable
(\cite{Polovinkin+Balashov,Rockafellar}).

By Corollary 2.\ref{S1} we have the estimate
$$
|(p_{1}-p_{2},x_{1}-x_{2})|\le \| p_{1}-p_{2}\|\|
x_{1}-x_{2}\|\le \Delta\ep(\Delta).
$$

From the equalities $(p_{i},x_{i})=s(p_{i},A)$, $i=1,2$, we
conclude that
$$
\ep_{1}(\| p_{1}-p_{2}\|)\| p_{1}-p_{2}\| = (p_{1},x_{1}-x_{2}),
\qquad \ep_{2}(\| p_{1}-p_{2}\|)\| p_{1}-p_{2}\| =
(p_{2},x_{2}-x_{1}).
$$
By the formulae
\begin{equation}\label{*10}
(p_{1},x_{1}-x_{2})+(p_{2},x_{2}-x_{1})=(p_{1}-p_{2},x_{1}-x_{2}),
\end{equation}
and $(p_{1},x_{1}-x_{2})\ge 0$, $(p_{2},x_{2}-x_{1})\ge 0$ we have
$$
\max\{ |(p_{1},x_{1}-x_{2})|, |(p_{2},x_{2}-x_{1})|\}\le
(p_{1}-p_{2},x_{1}-x_{2})\le \ep (\Delta)\Delta.
$$
\qed

It is well known that the external polyhedral approximation $\hat
A $ of the convex compact set $A\subset\R^{n}$ with supporting
function $s(p,A)$ on the grid $\C$ satisfies the formula
$s(p,\hat A) = \co \CC s(p,A)$ (\cite{sca,Polovinkin+Balashov}).

\begin{Th}\label{T2}
 Let
$A\subset\R^{n}$ be a convex compact set with the modulus of
convexity $\delta (\ep)$, $\ep\in [0,\diam A]$. Let $\C$ be a
grid with the step $\Delta\in (0,\frac12 )$, $\delta (\diam
A)/\diam A>\frac{\Delta}{4-\Delta^{2}}$. Then
$$
h(A,\hat A)\le \frac87 \ep (\Delta)\Delta,
$$
where $\ep (\Delta)$ is a solution of the equation $\frac{\delta
(\ep)}{\ep}=\frac{\Delta}{4-\Delta^{2}}$.
\end{Th}

\proof From the inclusion $A\subset \hat A$, formula $s(p,\hat
A)=\co\CC s(p,A)=\co\UU s(p,A)$ (see Lemma 2.\ref{L2}) and from
Definition 2.\ref{gridOp} it follows that (in terms of Definition
1.\ref{grid})
\begin{equation}\label{*12}
0\le s(p,\hat A)-s(p,A)\le \UU s(p,A)-s(p,A)= \alpha
\sum\limits_{i\in I_{p}}\hat\alpha_{i}(s(p_{i},A)-s(\hat p
,A)),\quad\forall p\in\R^{n}.
\end{equation}
Let $\hat x = \arg\max\limits_{x\in A}(\hat p,x)$.
$$
s(p_{i},A)-s(\hat p,A) = (\hat x, p_{i}-\hat p )+\ep_{i}(\|
p_{i}-\hat p\|)\| p_{i}-\hat p\|,
$$
and by properties of vector $\hat p$ (Lemma 2.\ref{L1}) and
Theorem 2.\ref{T1} we conclude that
$$
|\ep_{i}(\| p_{i}-\hat p\|)\| p_{i}-\hat p\| |\le
\ep(\Delta)\Delta,\quad \forall i\in I_{p}.
$$
Finally, using the equality $\sum\limits_{i\in
I_{p}}\hat\alpha_{i}p_{i}=\hat p$, we obtain that
$$
\sum\limits_{i\in I_{p}}\hat\alpha_{i}\left( s(p_{i},A)-s(\hat p
,A)\right) = \sum\limits_{i\in I_{p}}\hat\alpha_{i} \left( (\hat
x, p_{i}-\hat p )+\ep_{i}(\| p_{i}-\hat p\|)\| p_{i}-\hat p\|
\right)\le
$$
$$
\qquad\qquad\qquad\qquad\le\max\limits_{i\in I_{p}}| \ep_{i}(\|
p_{i}-\hat p\|)\| p_{i}-\hat p\||\le \ep (\Delta)\Delta.
$$
By formula (2.\ref{*12}) it follows
$$
0\le s(p,\hat A)-s(p,A)\le \alpha\ep (\Delta)\Delta = \frac{\| p
\|}{\|\hat p\|}\ep (\Delta)\Delta\le \frac{\ep
(\Delta)\Delta}{1-\frac12\Delta^{2}}\|p\|\le \frac87 \ep
(\Delta)\Delta\|p\|.
$$
By formula (1.\ref{*3}) we  get
$$
h(A,\hat A)\le \frac87 \ep (\Delta)\Delta.
$$
\qed

\begin{Corol}\label{Z1}
For any convex compact set $A$ the modulus of convexity $\delta$
satisfies the estimate $\delta (\ep)\le C\ep^{2}$. Thus the
typical value of $\ep (\Delta )$ is $\ep (\Delta
)\asymp\Delta^{s}$, $\Delta\to +0$, where $s\in (0,1]$.
\end{Corol}

\begin{Corol}\label{Z2}
The estimate of Theorem 2.\ref{T2} is exact.
\end{Corol}

Consider an example. Let $\R^{2}$ be the Euclidean plane with the
standard basis $Ox_{1}x_{2}$. Let $A=\{ x_{2}\ge
|x_{1}|^{s}\}\cap B_{1}(0)$, $s\ge 2$.

The modulus of convexity for the set $A$ equals $\delta (\ep) =
\ep^{s}/2^{s}$ for small $\ep>0$ and it is realized on the segment
$\left[ (-\frac{\ep}{2},\frac{\ep^{s}}{2^{s}}),
(\frac{\ep}{2},\frac{\ep^{s}}{2^{s}})\right]$. Let $a$ and $b$ be
two points from $\d A$:
$$
a=\left(-\frac{\ep}{2},\frac{\ep^{s}}{2^{s}}\right),\qquad
b=\left(\frac{\ep}{2},\frac{\ep^{s}}{2^{s}}\right).
$$
Let $p_{a}$ and $p_{b}$ be unit normals to the set $A$ at the
points $a$ and $b$, respectively. It is easy to calculate that
$$
p_{b}=\frac{\left(
s(\ep/2)^{s-1},-1\right)}{\sqrt{s^{2}(\ep/2)^{2(s-1)}+1}},\quad
p_{a}=\frac{\left(
-s(\ep/2)^{s-1},-1\right)}{\sqrt{s^{2}(\ep/2)^{2(s-1)}+1}},\quad
\| p_{b}-p_{a}\| =
\frac{2s(\ep/2)^{s-1}}{\sqrt{s^{2}(\ep/2)^{2(s-1)}+1}}.
$$
Let $\Delta = \| p_{b}-p_{a}\|$ and $p_{a}$, $p_{b}$ be adjacent
vectors of some grid $\C$ with the step $\Delta$. Suppose that
the grid $\C$ has symmetry with respect to the line $Ox_{2}$.

Then $\ep\asymp\Delta^{\frac{1}{s-1}}$ (for small $\Delta$). The
tangent line to the graph $x_{2}=|x_{1}|^{s}$ at the point $b$ is
$$
y_{tan}(x_{1})=s\left(\frac{\ep}{2}\right)^{s-1}\left(
x_{1}-\frac{\ep}{2}\right)+\left(\frac{\ep}{2}\right)^{s}.
$$
We have $y_{tan}(0)=-(s-1)\left(\frac{\ep}{2}\right)^{s}$. The
point $c=(0,y_{tan}(0))$ belongs to the set $\hat A$ (because
approximation $\hat A$ has symmetry with respect to the line
$Ox_2$). Hence
$$
h(A,\hat A)\ge \varrho (c,A) =
|y_{tan}(0)|=(s-1)\left(\frac{\ep}{2}\right)^{s}.
$$
So we have $h(A,\hat A)\ge C\cdot \Delta^{\frac{s}{s-1}}$. The
same order $\frac{s}{s-1}$ is given by Theorem 2.\ref{T2}.

\section{Approximation by presupporting functions}

Suppose that we know a presupporting function $f$ of a convex set
$A\subset\R^{n}$. We want to estimate the difference $\co \CC f
(p)- \co f (p)$ for all $p\in \R^{n}$. The last question is
equivalent to the question of evaluation of the value $h
(A,\tilde A)$ where $A=\{ x\ |\ (p,x)\le \co f(p),\ \forall
p\in\R^{n} \}$, $\tilde A=\{ x\ |\ (p,x)\le \co\CC f(p),\ \forall
p\in\R^{n} \}$.

The \it geometric difference \rm of sets $B,A\subset\R^{n}$ is the
set
$$
B\diff A = \{ x\ |\ x+A\subset B\}=\bigcap\limits_{a\in A}(B-a).
$$

We shall obtain the solution for two particular cases of
presupporting function. The first case is $f(p)=s(p,B)-s(p,A)$,
where $A$ and $B$ are convex compacta. In the case $B\diff
A\ne\emptyset$ the convex hull $\co f(p)$ equals the supporting
function $s(p,B\diff A)$ of the geometric difference $B\diff A$
\cite[Formula (1.11.18)]{Polovinkin+Balashov}.

The second case is $f(p)=\min\{ s(p,A),s(p,B)\}$, where $A$ and
$B$ are convex compacta. In the case $A\cap B\ne\emptyset$ the
convex hull $\co f(p)$ equals the supporting function $s(p,A\cap
B)$ of the intersection $A\cap B$ \cite[Formula
(1.11.17)]{Polovinkin+Balashov}.

The considered cases have important role for computational
geometry \cite{Dekstra} and for linear differential games
\cite{LDG,Pontryagin}.

\begin{Th}\label{T3}
Let $A,B\subset\R^{n}$ be convex compacta and suppose that $B$ is
uniformly convex with modulus $\delta_{B}$. Let
$f(p)=s(p,B)-s(p,A)$. Let $B_{r_{0}}(a)\subset B\diff A\subset
B_{d}(a)$. Let $\C$ be a grid with step $\Delta\in (0,\frac12 )$;
$\delta_{B} (\diam B)/\diam B>\frac{\Delta}{4-\Delta^{2}}$. Then
\begin{equation}\label{*13}
\co f(p)\le \co\CC f(p)\le \co
f(p)+\frac{8d}{7r_{0}}\ep_{B}(\Delta)\Delta\| p\|,\qquad\forall
p\in\R^{n},
\end{equation}
where $\ep=\ep_{B}(\Delta)$ is a solution of
$\frac{\delta_{B}(\ep)}{\ep}=\frac{\Delta}{4-\Delta^{2}}$.
\end{Th}

\proof The left inequality in (3.\ref{*13}) is obvious. By
\cite[Formula (1.11.18)]{Polovinkin+Balashov} we have
$$
\co f(p) = \co (s(p,B)-s(p,A))=s(p, B\diff A).
$$
Let $C=B\diff A$. By the formula
$$
\CC f(p)+\CC s(p,A) = \CC s(p,B)
$$
we obtain that
$$
\co\left( \CC f(p)+\CC s(p,A)\right)=\co\CC s(p,B)=s(p,\hat B),
$$
where $\hat B = \{ x\ |\ (p,x)\le \CC s(p,B),\ \forall
p\in\R^{n}\}$. Using inequality $\co (f+g)\ge \co f+\co g$ (which
is true for any functions $f$, $g$) we have
$$
\co\CC f(p)+\co\CC s(p,A)\le s(p,\hat B),
$$
or
$$
\co\CC f(p)\le s(p,\hat B)-\co\CC s(p,A)=s(p,\hat B)-s(p,\hat
A)\le s(p,\hat B)-s(p,A).
$$
By the last inequality
$$
\co\CC f(p)\le \co\left( s(p,\hat B)-s(p,A)\right)=s(p,\hat
B\diff A).
$$

Let $\| p\|=1$. Using (1.\ref{*3}) we have
$$
\co\CC f(p)-\co f(p)\le s(p,\hat B\diff A) - s(p,B\diff A)\le
h(\hat B\diff A,B\diff A).
$$
Let $h=h(B,\hat B)$. Using conditions of the theorem we conclude
that
$$
\begin{array}{l}
B\diff A\subset \hat B\diff A \subset \left(
B+B_{h}(0)\right)\diff A\subset \left(
B+\frac{h}{r_{0}}(C-a)\right)\diff A =\\ = \left(
B+\frac{h}{r_{0}}\left( B\diff A\right )\right)\diff
A-\frac{h}{r_{0}}a\subset \left( B+\frac{h}{r_{0}} B\diff
\frac{h}{r_{0}}A\right)\diff A-\frac{h}{r_{0}}a=\\
= \left( B\diff A\right)+\frac{h}{r_{0}}\left(\left( B\diff
A\right)-a\right)\subset  \left( B\diff A\right) +
\frac{h}{r_{0}}B_{d}(0).
\end{array}
$$
Hence
$$
h\left( B\diff A, \hat B\diff A\right)\le \frac{d}{r_{0}}h(B,\hat
B).
$$
Applying Theorem 2.\ref{T2} we finish the proof.\qed

\begin{Th}\label{T4}
Let $A,B\subset\R^{n}$ be uniformly convex compacta with moduli
$\delta_{A}$, $\delta_{B}$. Let $f(p)=\min\{s(p,A),s(p,B)\}$. Let
$B_{r_{0}}(a)\subset A\cap B$, $\max\{ \diam \hat A, \diam\hat
B\}\le d$. Let $\C$ be a grid with step $\Delta\in (0,\frac12 )$;
$\delta_{A} (\diam A)/\diam A>\frac{\Delta}{4-\Delta^{2}}$,
$\delta_{B} (\diam B)/\diam B>\frac{\Delta}{4-\Delta^{2}}$. Then
$\co f(p)\le$
\begin{equation}\label{*14}
\quad\le \co\CC f(p)\le \co f(p)+\frac{8}{7}\left( \max\{
\ep_{A}(\Delta), \ep_{B}(\Delta)\}+\frac{d}{r_{0}}\left(
\ep_{A}(\Delta)+\ep_{B}(\Delta)\right)\right)\Delta\| p\|, \
\forall p\in\R^{n},
\end{equation}
where $\ep=\ep_{X}(\Delta)$ is a solution of
$\frac{\delta_{X}(\ep)}{\ep}=\frac{\Delta}{4-\Delta^{2}}$, $X=A$
or $X=B$.
\end{Th}

\proof Let $C=A\cap B$ and $\hat C$ be external polyhedral
approximation of the set $C$ on the grid $\C$. By the inclusions
$C\subset A$, $C\subset B$ we have $\hat C\subset \hat A$, $\hat
C\subset \hat B$ and thus $\hat C\subset \hat A\cap\hat B$.

$$
\CC f(p) = \min\{ \CC s(p,A),\CC s(p,B)\},\qquad \co\CC f(p) =
\co \min\{ \CC s(p,A),\CC s(p,B)\}.
$$
By the inequality $\co\min\{ f,g\}\le \min\{ \co f,\co g\}$
(which is valid for any functions $f$, $g$) we have
$$
\co\CC f(p)\le \min\{ \co\CC s(p,A), \co\CC s(p,B)\} = \min\{
s(p,\hat A), s(p,\hat B)\}.
$$
Hence
$$
\co\CC f(p)\le \co\left( \min\{ s(p,\hat A), s(p,\hat B)\}\right)
=s\left(p, \hat A\cap\hat B\right).
$$
Let $\| p\|=1$. Using (1.\ref{*3}) we have
$$
\co\CC f(p)-\co f(p)\le s\left(p, \hat A\cap\hat B\right) -
s(p,A\cap B)\le h(A\cap B,\hat A\cap\hat B).
$$
Applying Theorem 3.1 \cite{Balashov+Repovs1} we obtain that
$$
h(A\cap B,\hat A\cap\hat B) \le \max\{ h(A,\hat A), h(B,\hat
B)\}+\frac{d}{r_{0}}\left( h(A,\hat A)+ h(B,\hat B)\right).
$$
Using of Theorem 2.\ref{T2} ends the proof.\qed

\section{On finding the convex hulls}

\begin{Th}\label{T5}
Let $A\subset\R^{n}$ be a uniformly convex compact set with
modulus of convexity $\delta$. Let
$$
r_{0}=\sup\{ r\ge 0\ |\ \exists\, a\in\R^{n}:\ B_{r}(a)\subset
A\}.
$$
Let a point $a\in\R^{n}$ be such that $B_{r_{0}}(a)\subset A$ and
$d=\sup\limits_{x\in A}\| x-a\|$. Then
$$
d\le \max\left\{2r_{0},
r_{0}+\delta^{-1}\left(\frac{r_{0}}{2}\right)\right\},
$$
where $\delta^{-1}$ is the inverse function for modulus $\delta$.
\end{Th}

\proof  Suppose that $d>2r_{0}$. We shall prove that $d\le
r_{0}+\delta^{-1}\left(\frac{r_{0}}{2}\right)$.

Let $b\in A$ and $\| a-b\|=d$. Let $L$ be any 2-dimensional
affine plane which contains points $a$, $b$. Let
$\L=(a-b)^{\perp}$, $\dim\L = n-1$.

Our further consideration will take place on the plane $L$. Let
the line $l$ be orthogonal to the line $\aff\{ a,b\}$ and $a\in
l$. Let $\{x,y\}=l\cap \d B_{r_{0}}(a)$. From the triangle $xab$
we have $\| x-b\|\ge d-r_{0}$, from the triangle $yab$ we have $\|
y-b\|\ge d-r_{0}$.

Let $z=\frac{a+b}{2}$ and let the line $l_{1}$ be parallel to the
line $l$ and $z\in l_{1}$. Let $x_{1}=l_{1}\cap \aff\{ x,b\}$,
$y_{1}=l_{1}\cap \aff\{ y,b\}$. By the uniform convexity of the
set $A$ we obtain that
$$
\left[ z,z+\frac{x_{1}-z}{\| x_{1}-z\|}\left(
\frac{r_{0}}{2}+\delta (\| x-b\|)\right)\right]\bigcup \left[
z,z+\frac{y_{1}-z}{\| y_{1}-z\|}\left( \frac{r_{0}}{2}+\delta (\|
y-b\|)\right)\right]\subset A,
$$
hence
$$
\left[ z,z+\frac{y_{1}-z}{\| y_{1}-z\|}\left(
\frac{r_{0}}{2}+\delta (d-r_{0})\right)\right]\bigcup \left[
z,z+\frac{x_{1}-z}{\| x_{1}-z\|}\left( \frac{r_{0}}{2}+\delta
(d-r_{0})\right)\right]\subset A.
$$
If $R=\frac{r_{0}}{2}+\delta (d-r_{0})>r_{0}$ then (due to the
previous inclusion being valid for any 2-dimensional plane $L$
with $\{ a,b\}\subset L$) we have
$$
B_{R}(z)\cap (\L+z)\subset A
$$
and thus
$$
\co\left( B_{r_{0}}(a)\cup \left(B_{R}(z)\cap
(\L+z)\right)\right)\subset A.
$$
By the last inclusion and by the inequality $\| a-z\|>r_{0}$ we
obtain that a shift of the ball $B_{r_{0}}(a)$ on a small distance
in the direction $b-a$ occurs in the interior of the set $A$.
Hence
 $r_{0}$ is not the maximal radius of balls from $A$. This contradiction
shows that $\frac{r_{0}}{2}+\delta (d-r_{0})\le r_{0}$.\qed



Now we describe an algorithm for finding the convex hull of a
positively uniform function \cite{Polovinkin+Balashov,LDG}.

Suppose that $\C$ is a grid with step $\Delta\in (0,\frac12)$,
$f(p)$ is a positively uniform continuous function and $\tilde A
= \{ x\ |\ (p,x)\le \CC f(p),\ \forall p\}$. We wish to calculate
$\co\CC f(p)$ for all $p\in \C$. In other words, we wish to find
$\CC \co\CC s(p,\tilde A)$. The problem can be solved as a
collection of problems of linear programming: for all $q\in \C$
to find
$$
(q,x)\to\max\quad (p,x)\le \CC f(p),\ \forall p\in \C.
$$

We shall describe an approximate algorithm from \cite{silin},
\cite[Theorem 2.6.3]{Polovinkin+Balashov}, \cite{LDG} and discuss
its error for the case of uniformly convex set $A=\{ x\ |\
(p,x)\le \co f(p),\ \forall p\in\R^{n}\}$ with modulus of
convexity $\delta$.

Suppose that $B_{r_{0}}(a)\subset A$ is the ball of maximum
radius in the set $A$ and $d=\sup\limits_{x\in A}\| x-a\|$.

We often do not know the precise values of $a$, $r_{0}$, but we
can easily calculate the ball of maximum radius $B_{R}(b)$,
$r_{0}\le R$, from $\tilde A$: it suffices to solve the following
problem of linear programming
$$
R\to\max\qquad (p,b)+R\le f(p),\quad\forall p\in\C.
$$
The solution $(b,R)\in\R^{n}\times\R$ gives the center of the
ball and its radius. In this case $B_{R}(b)\subset \tilde A$ and
$\tilde A \subset b+\left(d+\frac{4d^{2}}{r_{0}}\Delta\right)
B_{1}(0)$, see \cite[Corollary 2.6.2]{Polovinkin+Balashov}. Thus
\begin{equation}\label{rd}
(p,b)+R\| p\|\le f(p)\le (p,b)+\left(
d+\frac{4d^{2}}{r_{0}}\Delta\right)\| p\|,\qquad \forall
p\in\R^{n}.
\end{equation}

We shall further assume that $b=0$.

 The
first step of the approximate algorithm is to calculate for all
$q\in\C$ the values
$$
s^{\circ}(q)=\max\left\{ \frac{(p,q)}{f(p)}\ |\ p\in \C\right\}.
$$

The second step is to define $z(q)=q/s^{\circ}(q)$, $\forall
q\in\C$. Then the polyhedron
$$
A_{1}= \co\bigcup\limits_{q\in\C}z(q)
$$
is an approximation of $A$. The approximate value for $\co\CC
f(p)$, $p\in \C$, is
$$
\max\{ (p,z(q))\ |\ q\in \C\}.
$$
By \cite[Theorem 2.6.3]{Polovinkin+Balashov} we have under the
assumption (4.\ref{rd}) that $A_{1}\subset A$,
\begin{equation}\label{algEst}
h(A,A_{1})\le \frac{2\left(
d+\frac{4d^{2}}{r_{0}}\Delta\right)^{2}}{R}\Delta,
\end{equation}
and
$$
0\le \co\CC f(p) - \max\{ (p,z(q))\ |\ q\in \C\} \le \frac{2\left(
d+\frac{4d^{2}}{r_{0}}\Delta\right)^{2}}{R}\Delta,\qquad \forall
p\in\C.
$$

So the error of the algorithm is proportional to the step
$\Delta$ and to the value $\frac{1}{r_{0}}$ in the general case.

Consider the case when the set $A$ has modulus of convexity of
the second order: $\delta (\ep)=C\ep^{2}+o(\ep^{2})$, $\ep\to+0$.
Then under assumption $\Delta\in (0,\sqrt{r_{0}})$ we obtain by
Theorem 4.\ref{T5} that
$$
\frac{2\left( d+\frac{4d^{2}}{r_{0}}\Delta\right)^{2}}{R}\le
\frac{2\left(
d+\frac{4d^{2}}{r_{0}}\Delta\right)^{2}}{r_{0}}\le\mbox{Const}.
$$
and the error of the algorithm does not depend on the radius $R$
of an interior ball.

\section{Epilogue}

{\bf 1.} By Theorem 4.\ref{T5} we can estimate the value
$\frac{d}{r_{0}}$ in the theorems from Section 3. For example, in
Theorem 3.\ref{T3} in the case when $B_{r_{0}}(a)$ is the ball of
maximum radius from $B\diff A$ and $d=\sup\limits_{x\in B\diff
A}\| x-a \|$, we have (for small $r_{0}>0$)
$$
\frac{d}{r_{0}}\le\frac{r_{0}+\delta_{B}^{-1}\left(\frac{r_{0}}{2}\right)}{r_{0}}.
$$
If the modulus $\delta_{B}$ has the second order at zero then
$\frac{d}{r_{0}}\asymp\frac{1}{\sqrt{r_{0}}}$, $r_{0}\to+0$.

{\bf 2.} We want to point out that if the sets $A$, $B$ are
uniformly convex with moduli $\delta_{A}$, $\delta_{B}$
respectively and $\hat A$, $\hat B$ are polyhedral approximations
of $A$ and $B$ on a grid $\C$ with the step $\Delta$, then
$$
h(\hat A+\hat B, \widehat{A+B})\le \frac87\ep_{A+B}(\Delta)\Delta
$$
and in general in spaces of 3 or more dimensions $\hat A+\hat
B\subset  \widehat{A+B}$, but $\hat A+\hat B\ne  \widehat{A+B}$.
So the sum of approximations does not equal the approximation of
sum.

{\bf 3.} The results can easily be reformulated in any finite
dimensional Banach space. The only obstacle for the proofs is in
Lemma 2.\ref{L1} when we estimate $\| \hat p\|\ge
1-\frac12\Delta^{2}$. One must demand from the space and the grid
that
$$
C=\inf\limits_{\| p\|=1}\left\|\sum\limits_{i\in
I_{p}}\hat\alpha_{i}p_{i}\right\|\in (0,1).
$$
Then one must  replace denominator $4-\Delta^{2}$ by the new
$2+2C$ and coefficient $\frac87$ by the new $\frac{1}{C}$ in all
theorems.
\bigskip

 \centerline{Acknowledgements}
\bigskip

  This research was supported by SRA grants P1-0292-0101,
J1-2057-0101, and BI-RU/08-09/001. The first author was supported
by  RFBR grant 10-01-00139-a, ADAP project "Development of
scientific potential of higher school" 2.1.1/500 and project of
FAP "Kadry" 1.2.1 grant P938 and grant 16.740.11.0128. We thank
the referee for comments and suggestions.


\bigskip

\begin{thebibliography}{99}

\bibitem{Aubin} \au{J.-P. Aubin, I. Ekeland,} \tit{Applied Nonlinear
Analysis,} John Wiley \& Sons Inc. New York, 1984.

\bibitem{Balashov+Repovs1} \au{M. V. Balashov, D. Repov\v{s},} \tit{On the spliting problem for selections,}
J. Math. Anal. Appl. 355:1 (2009), 277-287.


\bibitem{Balashov+Repovs2} \au{M. V. Balashov, D. Repov\v{s},} \tit{Uniform convexity and the spliting problem for selections,}
J. Math. Anal. Appl. 360:1 (2009), 307-316.


\bibitem{Balashov+Repovs3} \au{M. V. Balashov, D. Repov\v{s},} \tit{Weakly convex sets and the modulus of nonconvexity,}
J. Math. Anal. Appl. 371:1 (2010), 113-127.

\bibitem{Dekstra} \au{E. W. Dijkstra,} \tit{A Discipline of Programming,} Prentice-Hall
Series in Automatic Computation, 1976.


\bibitem{Lindestrauss+tzafriri}
 \au{J. Lindenstrauss, L. Tzafriri,}
 \tit{Geometry of Banach Spaces - II. Functional Spaces,}
Springer-Verlag, Berlin, 1979.

\bibitem{silin} \au{G. B. Orlova, D. B. Silin,} \tit{Approximate calculus of the convex
hull for the positively uniform function,} Vestnik MGU, Ser. 15,
Computational Mathematics and Cybernetics, (1997), 2, 32 - 35 (in
Russian).

\bibitem{sca} \au{E. S. Polovinkin,} \tit{Strongly convex
analysis,}  Sbornik: Mathematics (1996), 187:2, 259-286.

\bibitem{Polovinkin+Balashov} \au{E. S. Polovinkin, M. V.
Balashov,} \tit{Elements of Convex and Strongly Convex Analysis,}
Fizmatlit, Moscow, 2007. (in Russian).

\bibitem{LDG} \au{E. S. Polovinkin, G. E. Ivanov, M. V. Balashov, R. V.
Konstantinov, A. V. Khorev,} \tit{An algorithm for the numerical
solution of linear differential games}. Sbornik: Mathematics.
2001, 192:10, 1515-1542.

\bibitem{Polyak} \au{B. T. Polyak,} \tit{Existence theorems and convergence of minimizing sequences
in extremum problems  with restrictions,} Soviet Math, 7 (1966),
72-75.

\bibitem{Pontryagin} \au{L. S. Pontryagin,} \tit{Linear differential games of pursuit}, Mat.
Sb. (N.S.) 112(154):3(7) (1980), 307-330.

\bibitem{Rockafellar} \au{R.T. Rockafellar,} \tit{Convex analysis.} Princeton University Press,
Princeton, NJ, 1970.

\end{thebibliography}
\end{document}